\documentclass[12pt]{article}
\usepackage{amsmath, amssymb, amsfonts}

\def\CC{{\rm \kern.24em \vrule width.02em height1.4ex
depth-.05ex \kern-.26em C}}

\def\TagOnRight

\def\AA{{\it I}\hskip-3pt{\tt A}}

\def\QQ{\rlap {\raise 0.4ex \hbox{$\scriptscriptstyle |$}}
  {\hskip -0.1em Q}}

\newcommand{\be}{\begin{equation}}
\newcommand{\ee}{\end{equation}}
\newcommand{\bea}{\begin{eqnarray}}
\newcommand{\eea}{\end{eqnarray}}
\newcommand{\Bea}{\begin{eqnarray*}}
\newcommand{\Eea}{\end{eqnarray*}}

\catcode`\@=11
\def\theequation{\@arabic{\c@equation}}
\catcode`\@=12

\newcommand{\bi}{\begin{itemize}}
\newcommand{\ei}{\end{itemize}}

\newtheorem{Definition}{Definition}[section]
\newtheorem{Theorem}[Definition]{Theorem}
\newtheorem{Lemma}[Definition]{Lemma}

\newtheorem{Corollary}[Definition]{Corollary}
\newtheorem{Remark}[Definition]{Remark}

\renewcommand{\theequation}{ \thesection.\arabic{equation} }

\begin{document}
\parskip 12pt
\parindent 0pt
\renewcommand{\thepage}{\arabic{page}}

\title{On the Feynman-Kac Formula
 }
\author{B.Rajeev\\ Email: brajeev@isibang.ac.in}
\maketitle
\begin{abstract} In this article, given $y :[0,\eta)\rightarrow H$ a
continuous map into  a Hilbert space $H$ we study the equation
\[\hat y(t) = e^{\int_0^tc(s,\hat y)}y(t)\] where $c(s,\cdot)$ is a given
`potential' on $C([0,\eta),H)$. Applying the transformation $y
\rightarrow \hat y$ to the solutions of the SPDE and PDE underlying
a diffusion, we study the Feynman-Kac formula.

\end{abstract}

\section{Introduction} One of the well known formulas at the
boundary of probability and analysis is the Feynman-Kac formula
$u(t,x) = E_x(f(X_t)e^{-\int_0^tV(X_s)ds})$, which represents the
solution $u(t,x)$ of the evolution equation for the operator $L - V$
where $L$ is the infinitesimal generator of a diffusion $(X_t, P_x),
x \in \mathbb R^d$,$V(x) \geq 0$ the potential function and $f$ the
initial value(\cite{MK}). We refer to \cite{BO1}, \cite{KS},
\cite{OK} for basic material on this topic. It is also known that
this formula defines a sub-Markovian semi-group whose underlying
process $(\hat X_t)$ is obtained from $(X_t)$ by the operation known
as `killing' according to the multiplicative functional $M_t :=
e^{-\int_0^tV(X_s)ds}$ (\cite{RY}). It maybe of interest therefore
to have an answer to the following natural question : is it possible
to have a `pathwise' construction of the process $\hat X$. The
special case when $(X_t)$ satisfies an It$\hat{o}$ stochastic
differential equation (SDE), is of interest. However, it turns out
that it is the SPDE satisfied by the distribution valued process
$(\delta_{X_t})$ (\cite{BR1}) rather than the SDE for $(X_t)$ that
is more relevant for our purposes.

 To motivate our 'pathwise' construction, we proceed as follows.
 Let $H$ be a separable real Hilbert space
and consider $C([0,\eta),H)$, the space of continuous functions on
$[0,\eta), 0 < \eta \leq \infty$, with values in $H$. Let $u(t) \in
C([0,\eta),H)$ be the solution of the following evolution equation
in $H$ viz.
\[\partial_t u(t) = Lu(t)\] with $L : H \rightarrow H$ say, a linear
operator. Consider $\bar u(t) := u(t)e^{\int_0^tc(s,u)ds}$, where
$c(s,\cdot): C([0,\eta),H) \rightarrow \mathbb R$ is a given
function (the potential). Then, integrating by parts, it is easy to
see that $\bar u$ solves
\[\partial \bar u(t) = L\bar u(t) + c(t,u)\bar u(t). \] We would
have a good and proper evolution equation for $\bar u(t)$ if we were
able to write $c(t,u) = \hat c(t,\bar u)$. If the map $S(u) := \bar
u = u(t)e^{\int_0^tc(s,u)}ds$ were invertible, then we may define $
\hat c(t,u) := c(t, R(u))$ where $R(u):= S^{-1}(u)$ so that $\hat
c(t, \bar u) = c(t,u).$ It is easy to see that the inverse $R$ is a
path transformation $R: C([0,\eta),H) \rightarrow C([0,\eta),H)$
induced by the `potential' $c: [0,\eta)\times C([0,\eta),H)
\rightarrow \mathbb R$ as follows : For a given $y \in
C([0,\eta),H), R(y) \in C([0,\eta),H)$ is the solution $\hat y$ of
the equation
$$ \hat y(t)= y(t) e^{- \int\limits_0^t c(s,\hat y)ds}.
$$  In section 2, we prove existence and uniqueness to the above
equation in Theorem (2.2), using a fixed point argument. Thus the
map $R$ is well defined and injective. Since $-c$ satisfies the
conditions of Theorem (2.2) whenever $c$ does, the map $R$ is also
onto. From a modeling point of view, $R(y)$ maybe viewed as a
perturbation, induced by the potential $c(t,y)$, of the trajectory
of a particle represented by $y(.)$. We deal with real Hilbert
spaces as we consider applications only to the theory of diffusions.
However, complex Hilbert spaces and complex valued potentials (with
the corresponding interpretation of `amplitude' and  `phase') are
also of interest.

Given a diffusion $(X_t,0 \leq t < \eta,P_x, x \in \mathbb R^d)$, we
try to realise the Feynman-Kac formula by applying the above
transformation to the paths of the diffusion.  We remark here that
we could choose $H = \mathbb R^d$ but this does not lead to the
Feynman-Kac formula (see Remark (3.2)). However, if we look at the
process $ (Y_t) := (\delta_{X_t^x})$, upto time $\eta$, then this is
a semi-martingale in a Hilbert space ${\cal S}_p$- the so called
Hermite-Sobolev space - and indeed is the unique solution of a quasi
linear stochastic partial differential equations (SPDE)
(\cite{BR1},\cite{BR2}); one may then look at the process $(\hat
Y_t):= (e^{-\int_0^tV(X_s)ds}\delta_{X_t})$ and using the rules of
stochastic calculus write an SPDE for $\hat Y$. Note that we can
write $V(X_s) = \langle V,\delta_{X_s}\rangle =:
c(s,\delta_{X_\cdot}), 0 \leq s < \eta $, if $V$ belongs to a
suitable class of test functions.

In section 3, we show that when $(Y_t)$ satisfies a quasi linear
SPDE in ${\cal S}_p$ then  $\hat Y_t$ is the solution a new SPDE
with a potential term viz. $c(t,\hat Y)$ and whose coefficients are
defined on the path space $C([0,\eta),{\cal S}_p)$ using the
coefficients of the original equation and the transformation
discussed above. This transformation works at both levels viz. the
SPDE and the PDE underlying the diffusion, although the `Kac
functional' (we use the terminology from \cite{CV}) induced by the
potential function $V(x)$ is necessarily different in the two cases
(see the discussion on diffusions in section 5). In section 4, we
allow $ c(.,.)$ to depend also on $x \in \mathbb R^d$ and we show
that the above transformation may also be applied directly to the
solutions of a class of non-linear PDE's. We conclude in Section 5
with a discussion on two classes of examples in both of which the
functional $c(t,x,y)$ depends on $x$ albeit in different ways. The
second example that we discuss in Section 5 concerns diffusion
processes and shows also the connections that can arise between the
transformations of the solutions to the SPDE and the associated PDE.
In Sections 3,4 and 5, we work in the framework of \cite{BR2} to
which we refer for results relating to SPDE's, the related notations
and references. See also Example 7 of \cite{BR2} where we had
briefly indicated the results in Section 2.

\section{A Transformation on path space}
{ \rm Let $H$ be a separable real Hilbert space with norm denoted by
$\|.\|$. We consider for $0 \leq T < \infty$, the space $C([0,T],
H)$ of continuous functions $y : [0, T] \rightarrow H$ with the
sigma field ${\cal B}_t,  0 \leq t \leq T$ generated by the
coordinate maps upto time $t$. For a continuous map $y: [0,T]
\rightarrow H$ and $0 \leq s \leq T$, we denote its norm on
$C([0,s],H)$ by $\|y\|_s :=\sup\limits_{u\leq s} \|y(u) \|$. Fix
$T>0$. Let $c : [0,T]\times C([0,T],H) \rightarrow \mathbb{R}$
satisfy
\begin{enumerate}
\item{} For $0 \leq t \leq T$, ~~~$|c(t, y_1)-c(t,y_2)|~ \leq ~\beta \|y_1 -y_2\|_t$ \\
where $\beta=\beta(T)$ depends only on $T$.

We note that as a consequence of the condition 1) we have the
following : for $0 \leq s \leq T$ and $y_1, y_2 \in C([0,T],H),~
y_1(u) = y_2(u), 0 \leq u \leq s$ implies $c(s,y_1) = c(s,y_2)$.
\item{} For $\alpha > 0$ and $T > 0$ there exists a constant $M(\alpha ,
T)$ such that
\[|c (t,y)| \leq M(\alpha,T) \]
for $0 \leq t \leq T$ and for all $y \in B(0,\alpha) \equiv
B(0,\alpha,T):= \\ \{y \in C([0,T],H), \|y\|_T \leq \alpha \} $.
\end{enumerate}

 We note that if $c(t,y)$ satisfies conditions 1 and 2 then
so does $-c(t,y)$.

Let $\alpha (t) \equiv \alpha (t,y) := e^{-\int\limits_0^t
c(s,y)ds}$ for $y \in C([0,T],H)$. Given a $\hat y \in
C([0,\eta),H)$ for some $\eta
>0$, and  $0 < T < \eta$ we consider the following equation in $C([0,T],H)$,
viz.
\begin{eqnarray}
\hat y(t)=  y(t)\alpha (t,\hat y)= y(t) e^{- \int\limits_0^t
c(s,\hat y)ds}
\end{eqnarray}
for $0\leq t \leq T$. We first derive an apriori estimate for the
distance between two solutions corresponding to two `inputs' $y_1$
and $y_2$.

\begin{Lemma} Let $y_1, y_2 \in C([0,\eta),H)$ and suppose $\hat y_1, \hat y_2$
are the corresponding solutions of (2.1). Then for every $0 < T <
\eta$,  we have the following estimate viz.
\begin{eqnarray}  \|\hat y_1 - \hat y_2\|_T
\leq M\|y_1 - y_2\|_T e^{M\|y_2\|_T\beta e^{\delta}}.\end{eqnarray}
where $\delta
> \beta T \|\hat y_1- \hat y_2\|_T $ and $M := e^{\int_0^T|c(s,\hat y_1)|}ds$.
\end{Lemma}

{\bf Proof} Let $0 < T < \eta $ and $\delta ,M$ as above. Then
\[ \int_0^T |c(s,\hat y_1)-c(s,\hat y_2)| ds \leq \beta T\|\hat y_1
- \hat y_2\|_T < \delta \] and consequently, using the elementary
estimate $|1-e^x| \leq e^{\delta}|x|, |x| < \delta$ we have for any
$0 \leq t \leq T$,\[ |(1- e^{\int_0^t c(s,\hat y_1)-c(s,\hat
y_2)ds})| \leq e^{\delta}\beta \int_0^T \|\hat y_1 - \hat y_2\|_s
ds.\] Then we have

\begin{eqnarray*}
\| \hat y_1-\hat y_2\|_T &=& \|(y_1-y_2)e^{-\int_0^{\cdot}c(s,\hat
y_1)ds} + y_2 e^{-\int_0^{\cdot}c(s,\hat y_1)ds}(1-
e^{\int_0^{\cdot} c(s,\hat y_1)-c(s,\hat y_2)ds})\|_T \\
&\leq & \|(y_1-y_2)\|_T e^{\int_0^{T}|c(s,\hat y_1)|ds} + \|y_2\|_T
e^{\int_0^{T}|c(s,\hat y_1)|ds} \\ &&\times \sup\limits_{t \leq
T}|(1- e^{\int_0^{t} c(s,\hat y_1)-c(s,\hat y_2)ds})| \\ &\leq &
M\|y_1-y_2\|_T + M\|y_2\|_T e^{\delta}\beta \int_0^T \|\hat y_1 -
\hat y_2\|_s ds \\ &\leq& M \|y_1-y_2\|_T e^{TM\|y_2\|_T\beta
e^{\delta}}
\end{eqnarray*} where the last step follows from Gronwal's
inequality. $\hfill{\Box}$

Let $ y_1 \in C([0,t_1],H), y_2 \in C([0,t_2],H)$ where $0 \leq t_1
< t_1+ t_2 < T$. In the proof of the following theorem we need the
following construction of `concatenation' $y_1 \diamond y_2 \in
C([0,T],H)$ of the paths $y_1 $ and $y_2$ : \Bea
 y_1\diamond y_2(s) &:=& y_1(s)|_{[0,t_1]}(s)+ (y_2(s-t_1)-y_2(0)
 + y_1(t_1))|_{(t_1,t_1+t_2]}(s) \\
&+& (y_2(S)-y_2(0) + y_1(t_1))|_{(t_1+t_2,T]}(s), \Eea } where $|_A$
is the indicator of the set $A$. The following theorem is our main
result.

\begin{Theorem} Let $\eta > 0$ and
let $c(t,y)$ satisfy conditions 1 and 2 above for every $T, 0 \leq T
< \eta$. Then for a given $y \in C([0,\eta),H)$ there exists a
unique $\hat y \in C([0,\eta), H)$ satisfying equation (2.1) for
every $T, 0 < T < \eta$.
\end{Theorem}

{\bf Proof:} {\rm It suffices to show existence and uniqueness of
equation (2.1) on $[0,T]$ for every $T < \eta$. Using uniqueness,we
can then patch up the solutions on overlapping intervals to get the
required solution. So let $0 < T < \eta$. Uniqueness is immediate
from (2.2).

To show existence on $[0,T]$, suppose $(0,T] =
\bigcup\limits_{n=0}^{m-1} (T_n,T_{n+1}]$. Fix $n, 0 \leq n \leq
m-1$. Define $\hat y(0) := y(0).$ Suppose that $\hat y(t), t \in
[0,T_n]$ has been defined. Then we extend $\hat y$ to the interval
$(T_n,T_{n+1}]$ as follows : We first solve the following equation
on $[0,T_{n+1}-T_{n}]$ viz. \bea \hat y_n(t) &=& y(t+T_n)\alpha(T_n,
\hat y)\alpha_n(t,\hat y_n) - y(T_n)\alpha(T_n,\hat y) \\   &=&
y_n(t)\alpha_n(t,\hat y_n) + a_n \nonumber\eea where for $ y \in
C([0,T_{n+1}-T_n],H) {\rm ~ and ~} t \in [0,T_{n+1}-T_n],$
\[\alpha_n (t, y ) := e^{-
\int\limits_0^t c(s + T_n,\hat y\diamond y)ds},~~ y_n (t):=
y(T_{n}+t)\alpha (T_n, \hat y),\] and $a_n := -
y(T_n)\alpha(T_n,\hat y).$ \\ We extend $\hat y$ as follows :
\[\hat y(t) := \hat y_n(t-T_n)+ \hat y(T_n),~~ t \in (T_n,
T_{n+1}].\] Then provided $\hat y$ satisfies equation (2.1) in
$[0,T_n]$, we have for $t \in (T_n, T_{n+1}]$, \Bea \hat y(t) &:=&
\hat y_n(t-T_n)+ \hat y(T_n)
\\ &=& y(t)\alpha(T_n, \hat y)\alpha_n(t-T_n,\hat y_n) -
y(T_n)\alpha(T_n,\hat y) + \hat y(T_n)\\ &=&  y(t)\alpha(T_n, \hat
y)\alpha_n(t-T_n,\hat y_n)
\\ &=&  y(t)\alpha(t, \hat y)\Eea where in the third equality
we have used the assumption that $\hat y$ satisfies equation (2.1)
in $[0,T_n]$. As for the fourth equality, we use the fact that $\hat
y$ on the interval $(T_n,T_{n+1}]$ is the concatenation of $\hat y$
in $C([0,T_n],H)$ and $\hat y_n$ in $C([0,T_{n+1}-T_n],H)$ i.e.
$\hat y(t) = \hat y \diamond \hat y_n (t), t \in (T_{n}, T_{n+1}]$.

Thus it suffices to solve (2.3) on $[0,T_{n+1}-T_n]$ for a suitable
choice of $\{0= T_0 < T_1 < \cdots < T_n < \cdots < T_m = T\}$.

Let $\alpha
> \sup\limits_{s \leq T}\|y(s)\|,$ and $c(.,.)$ satisfy
conditions 1 and 2 on $[0,T]$ for some $M(\alpha):= M(\alpha,T)
{\rm~ and ~}\beta$. Let $\epsilon
> 0$ be such that $\epsilon e^{M(3\alpha)T} < \frac{\alpha}{2}$.
By uniform continuity of $y$ on $[0,T]$ we can divide $[0,T]$ into a
finite number (say $m$) of subintervals $[T_n,T_{n+1}]$, with $T_m =
T$ such that
\[\|y(t_1)- y(t_2)\| ~\leq~ \epsilon ~~~~ \forall t_1,t_2 \in [T_n,T_{n+1}],~ n = 0,\cdots m-1.\]

Next we choose $\delta > 0$ such that $|e^x - 1| < e^{\delta}|x|$
for $|x| < \delta$.

By refining the partition if necessary we may assume without loss of
generality that \[\alpha
M(3\alpha)e^{M(\alpha)T+\delta}(T_{n+1}-T_n) < \frac{\alpha}{2}; \]
\[K_n:= 2\alpha \beta e^{M(3\alpha)T+\delta}(T_{n+1}-T_n) < 1,~~ n =
0,\cdots, m-1 \] and,
\[ 2M(3\alpha))(T_{n+1}-T_n) < \delta.\]

With this choice of the partition $\{T_n\}$ we now solve equation
(2.3) on $[0,T_{n+1}-T_n]$ by a fixed point argument. Let $\alpha$
be as above. Recall the definition of $B(0,\alpha)$ from condition 2
above, with $T$ there replaced with $T_{n+1}-T_n$. For $z \in
B(0,\alpha), t \in [0,T_{n+1}-T_n]$ let
\[S_n(z)(t):= y(t+T_n)\alpha(T_n, \hat y)\alpha_n(t,z) -
y(T_n)\alpha(T_n,\hat y).\] Note that $\alpha_n(t,z)$ depends on
$\hat y \diamond z$ where $\hat y$ is the solution of (2.1) on
$[0,T_n]$. Assume that $\hat y \in B(0,\alpha)$. Then we claim that
\[ S_n : B(0,\alpha) \subset C([0,T_{n+1}-T_n],H)\rightarrow
B(0,\alpha).\] To see this we write $S_n(z)(t)$ as \[ S_n(z)(t) =
(y(t+T_n) - y(T_n))\alpha(T_n, \hat y)\alpha_n(t,z)
+y(T_n)\alpha(T_n,\hat y)(\alpha_n(t,z) -1).\] Let $t \in
[0,T_{n+1}-T_n]$. Then from the triangle inequality and the choice
of $\epsilon$ and $\{T_n\}$ we have

\Bea \|S_n(z)(t)\| &\leq & \|(y(t+T_n) - y(T_n))\| \alpha(T_n, \hat
y)\alpha_n(t,z)
\\ && + \|y(T_n)\|\alpha(T_n,\hat y)|(\alpha_n(t,z) -1)| \\
&\leq& e^{M(3\alpha)T}\|(y(t+T_n) - y(T_n))\|\\ &+& \alpha
e^{M(\alpha)T+\delta}|\int_0^t c(u+T_n, \hat y \diamond z)du|
\\ &\leq & \epsilon e^{M(3\alpha)T} + \alpha M(3\alpha)
e^{M(\alpha)T+\delta}(T_{n+1}-T_n) \leq \frac{\alpha}{2} +
\frac{\alpha}{2} = \alpha. \Eea Note that in the second equality we
have used the fact that condition 2) implies $|\int_0^t c(s+T_n,
\hat y\diamond z)ds| < M(3\alpha)(T_{n+1}-T_n) < \delta , t \in
[0,T_{n+1}-T_n]$ and in the second and third inequality above we
have used the fact that $\|\hat y\diamond z(t)\| \leq 3\alpha , t
\in [0,T]$.

We now show that the map $S_n: B(0,\alpha) \rightarrow B(0,\alpha)$
is a contraction. Let $y_1,y_2 \in B(0,\alpha)$ and $S_n(\cdot)$ be
as defined above. For $t \in [0,T_{n+1}-T_n]$,
\begin{eqnarray*}
\|S_n(y_1) (t)-S_n(y_2)(t)\| &=& \|y (t+T_n)\| \alpha (T_n,\hat
y)\alpha_n(t,y_2)\\ && \times \left|
e^{\int\limits_0^t (c(s+T_n,\hat y \diamond y_1)-c(s+T_n,\hat y \diamond y_2))ds}-1\right|\\
&\leq& \|y(t+T_n) \| \alpha(T_n,\hat y)\alpha_n(t,y_2)
\\ && \times e^{\delta}\left|
\int\limits_0^t (c(s+T_n,\hat y \diamond y_1)-c(s+T_n,\hat y \diamond y_2))ds \right|\\
&\leq& e^{TM(3\alpha)+\delta} \alpha \beta 2(T_{n+1}-T_n)\|y_1 - y_2\|_{T_{n+1}-T_n}\\
\end{eqnarray*} and by definition of the constant $K_n$ we have
\[\|S (y_1)-S(y_2)\|_{T_{n+1}-T_n} \leq K_n\| y_1 -
y_2\|_{T_{n+1}-T_n}.\] Since $K_n < 1$ by our choice, the map \[y
\rightarrow S_n(y) :C([0,T_{n+1}-T_n],B(0,\alpha)) \rightarrow
C([0,T_{n+1}-T_n],B(0,\alpha))\] is a contraction on a complete
metric space and has a unique fixed point. Thus equation (2.3) has a
unique solution. This completes the proof of the Theorem.
$\hfill{\Box}$

\begin{Corollary} For $y \in C([0,\eta),H)$, let $R(y):= \hat y$
where $\hat y$ is the solution of (2.1). Then $R$ is one to one and
onto. Further for every $t > 0$, $R : C([0,t],H) \rightarrow
C([0,t],H)$ is a homeomorphism. In particular, for every $t > 0$,
the map $R : (C([0,t],H),{\cal B}_t) \rightarrow (C([0,t],H),{\cal
B}_t)$ is a measurable isomorphism.
\end{Corollary}
{\bf Proof}: To see that $R$ is one-one, suppose that $R(y_1)=
R(y_2)$. Then since this implies $\hat y_1 = \hat y_2$, we also have
$y_1 = y_2$. That $R$ is onto follows from the observation that if
$\hat y \in C([0,\eta),H)$ is given and if we define $y(t) := \hat
y(t) e^{\int_0^t c(s,\hat y)ds} $ then clearly $R(y) = \hat y$.

 Note that
for a given $y \in C([0,\eta),H)$, $R^{-1}(y) = \bar y :=
ye^{\int_0^{\cdot}c(s,y)ds}$ follows since $R(\bar y) = y$. Since
$R^{-1}$ has the same form as $R$ it suffices to show that $R$ is
continuous. But this is clear from (2.2). The last statement follows
from the continuity of $R$ and the fact that the Borel sigma field
on $C([0,t],H)$ is the same as ${\cal B}_t$. $\hfill{\Box}$

\section{Application to Stochastic PDE's}
Let ${\cal S}_p, p \in \mathbb R$ be the family of Hermite-Sobolev
spaces; $\cal S, \cal S'$ respectively the Schwartz space of rapidly
decreasing smooth functions and its dual. We refer to \cite{BR2},
\cite{KI2},\cite{KX} for the results and notations related to these
spaces that we use. We refer to \cite{DZ} and \cite{GM} for results
on stochastic calculus in Hilbert spaces. We work on a probability
space $(\Omega, {\cal F},P)$ on which is given an $r$-dimensional
standard Brownian motion $(B_t)$. Let $({\cal F}_t^B)_{t \geq 0}$ be
the filtration of $(B_t)$.
 We now consider solutions of the SPDE
\begin{eqnarray}
dY_t &=&  L (t,Y)dt +  A(t, Y) \cdot dB_t \\
Y_0 &=& Y. \nonumber
\end{eqnarray}
where $L, A_i, i = 1, \cdots r$ are second order quasi-linear
partial differential operators with coefficients $\sigma_{ij}, b_i :
{\cal S}_p \rightarrow \mathbb R^d, i = 1, \cdots d, j = 1, \cdots
r$ defined as follows
\begin{eqnarray*}
L(y) &:=& \frac{1}{2} \sum\limits^d_{i,j=1} ( \sigma
\sigma^t)_{ij}(y) \partial^2_{ij} y
 -\sum\limits^d_{i=1}  b_i (y) \partial_i y \\
 A_i(y) &:=& -\sum\limits^d_{j=1}\sigma_{ji}(y)
\partial_jy
\end{eqnarray*}

In \cite{BR2} we have proved existence and uniqueness of solutions
to equation (3.4) and shown that for a given $Y : \Omega \rightarrow
{\cal S}_p$ a unique solution ($Y_t,\eta$) exists under a Lipschitz
condition on the coefficients $\sigma_{ij}$ and $b_i$. Here $0 <\eta
\leq \infty$ is the lifetime of the process and if $\sigma_{ij},b_i$
are uniformly bounded on ${\cal S}_p$ then $\eta = \infty$ almost
surely (see \cite{BR2}, Proposition (5.2)).

Let $c(.,.) : [0,\infty)\times C([0,\infty),{\cal S}_p) \rightarrow
\mathbb{R}$ satisfy conditions 1 and 2 of Section 2 on bounded
intervals $[0,T]$. Let $\eta > 0$. Given $y \in C([0,\eta),{\cal
S}_p) $ let $\hat y$ be the solution of equation (2.1) given by
theorem (2.1) with $H = {\cal S}_p$ .

Suppose now that we are given $\sigma_{ij}, b_i :{\cal S}_p
\rightarrow \mathbb{R}$ and $L,A^i,i=1,\cdots r$ as above. The
transformation $y \rightarrow \hat y $ induced by the map $c(.,.)$
and equation (2.1) induces a corresponding transformation of maps
$\sigma_{ij}(.),b_i(.) \rightarrow \hat \sigma_{ij}(.,.), \hat
b_i(.,.)$ as follows : $\hat \sigma_{ij}, \hat b_i : [0,\eta) \times
C([0,\eta),{\cal S}_p) \rightarrow \mathbb{R}$ by $\hat \sigma_{ij}
(s,y) :=\sigma_{ij} (\hat y(s)),~~\hat b_i (s,y):=b_i(\hat y(s))$.
Define $\hat c(s,y):= c(s, \hat y),~ 0 \leq s < \eta,~ y \in
C([0,\eta),{\cal S}_p)$. Let $\hat L(t,y)$ and $\hat A_i(t,y)$ be
maps from $[0,\eta)\times C([0,\eta),{\cal S}_p) $ to ${\cal S}_p$
for fixed $\eta
>0$ defined as follows:
\begin{eqnarray*}
\hat L(s,y) &:=& \frac{1}{2} \sum\limits^d_{i,j=1} (\hat \sigma \hat
\sigma^t)_{ij}(s,y) \partial^2_{ij} y_s
 -\sum\limits^d_{i=1} \hat b_i (s,y) \partial_i y_s + \hat c(s,y)y_s\\
\hat A_i(s,y) &:=& -\sum\limits^d_{j=1} \hat \sigma_{ji} (s,y)
\partial_j y_s
\end{eqnarray*}

 Let $(Y_t,\eta)$ be a pathwise
unique strong solution of equation (3.4) with initial value $Y$.
Then for each $\omega \in \Omega$, the trajectory $Y_\cdot(\omega)
\in C([0,\eta(\omega)),{\cal S}_p)$. Define for $0 \leq t <
\eta(\omega)$
\[
\hat Y_t(\omega):= Y_t(\omega)e^{\int\limits_0^t
c(s,Y(\omega))ds}.\] Let $\hat \sigma_{ij},\hat b_i,\hat c, \hat L,
\hat A_i$ be as above. We take $\hat Y_t(\omega):= \delta , t \geq
\eta$ where $\delta$ is the coffin state. By the continuity of
$c(.,.)$ and the definition of a strong solution (see \cite{BR2},
$(\hat Y_t)$ is a continuous ${\cal F}_t^B$-adapted, $\hat {\cal
S}_p := {\cal S}_p \bigcup \{\delta\}$ valued process.}
\begin{Theorem} Let $(Y_t)_{0 \leq t < \eta}$ be a strong solution
of equation (3.4) and let $c(.,.)$ satisfy conditions 1 and 2 of
Section 2. Then $(\hat Y_t)_{0 \leq t < \eta}$ is a strong solution
of the equation
\begin{eqnarray}
d\hat Y_t &=& \hat L (t,\hat Y)~dt +\hat A(t,\hat Y) \cdot dB_t \\
\hat Y_0 &=& Y. \nonumber
\end{eqnarray}
If equation (3.4) has a unique strong solution, then so has equation
(3.5).
\end{Theorem}

{\bf Proof:}{\rm~ Let $M_t := e^{\int\limits_0^t c(s, Y)ds}$. To
prove existence, we use integration by parts. Indeed one can verify
the following equation by acting on it with a test function. We have
in differential form
\begin{eqnarray*}
d\hat Y_t &=& d(M_t Y_t) =Y_t~~ dM_t +M_t dY_t\\
&=&  \hat Y_t c(t,Y)~dt +L(Y_t) M_t~ dt +M_t A(Y_t)\cdot dB_t
\end{eqnarray*}
Now from the definition of $\hat Y_t{(\omega)}$, we have that for
each fixed $\omega$, $Y_t{(\omega)}, 0 \leq t < \eta(\omega)$ is the
unique solution $\hat y$ of equation (2.1) with $y(t) := \hat Y(t),
0 \leq t < \eta$, viz.
\[
\hat y(t) = \hat Y_t(\omega)e^{- \int\limits_0^t c(s,\hat y)ds}.
\]
It follows that $\sigma_{ij} (Y_t)=\hat \sigma_{ij} (t,\hat Y),
c(t,Y) = \hat c(t,\hat Y)$ etc. and hence from above,
\[
d\hat Y_t =\hat L(t,\hat Y)~ dt +\hat A(t,\hat Y) \cdot dB_t.
\]
The uniqueness of solutions of equation (3.5) follows from the
uniqueness of equation (3.4). Indeed if $\hat Y^1, \hat Y^2$ are
solutions of equation (3.5), then if $(Y^i_t)$ solves $Y^i(t) = \hat
Y^ie^{- \int\limits_0^t c(s,Y^i)ds}$ it is easy to check using the
integration by parts formula for the product $\hat Y^ie^{-
\int\limits_0^t c(s,Y^i)ds}$ and the definition of the `hat'
functionals that $Y^i, i=1,2$ both solve equation (3.4) and hence
$Y^1 =Y^2$ which in turn implies $\hat Y^1 =\hat Y^2$.
$\hfill{\Box}$

\begin{Remark}
Let $(X_t^x, \eta)$ be the solution of an $It\hat{o}$ SDE with
diffusion and drift coefficients $\bar \sigma_{ij}$ and $\bar b_i$
respectively and initial value $x \in \mathbb R^d$. Let  $\bar V :
\mathbb R^d \rightarrow \mathbb R$ be a locally bounded function. We
can apply Theorem 2.2, with $H = \mathbb R^d$ to transform $X$ into
$\hat X := T(X)$ with $T : C([0,\eta),\mathbb R^d) \rightarrow
C([0,\eta),\mathbb R^d)$ given by $T(X)(s):= e^{ \int\limits_0^t
c(s,X)ds}X_s, 0 \leq t < \eta$ with $c: [0,\infty)\times
C([0,\infty),\mathbb R^d)\rightarrow \mathbb R , c(s,y) := \bar
V(y_s)$ and $y \in C([0,\infty),\mathbb R^d)$. Then $\hat X$ will
satisfy an SDE with path dependent coefficients which can be
determined as in the case of $\hat Y$ in Theorem (3.1). In general,
the transformation $T$ applied to an $It\hat{o}$ process $(X_t)$
changes the drift term by adding a term like $c(t,\hat X)\hat X_t$.

On the other hand,let $Y_t := \delta_{X_t^x}, t < \eta$ with $Y_0 =
\delta_x$ and suppose that the coefficients $\bar \sigma_{ij}, \bar
b_i, \bar V \in {\cal S}_p, p > \frac{d}{4}$ . Let $\sigma_{ij},
b_i, V$ be the linear functionals on ${\cal S}_p$ given by $
\sigma_{ij}(y) := \langle \bar \sigma_{ij}, y\rangle $ etc. Then
$(Y_t)$ is the unique solution of the  SPDE (3.4) with $Y_0 =
\delta_x$ and with $c: [0,\infty)\times C([0,\infty),{\cal
S}_{-p})\rightarrow \mathbb R$ given by $c(s,y) := \langle \bar V,
y_s\rangle ,y \in C([0,\infty),{\cal S}_{-p})$, $\hat Y$ is the
unique solution of (3.5) upto the lifetime $\eta$.
\end{Remark}

\section{Application to PDE's} We now apply the transformation $y \rightarrow
\hat y$ developed in Section 2, to solutions of partial differential
equations of the form \begin{eqnarray}
\partial_t u(t,x) &=& L(x,u(t,x))\\
u(0,x) &=& u(x) \nonumber .
\end{eqnarray}
with $u : \mathbb R^d \rightarrow {\cal S}_p$. Here the operator
$L(x): {\cal S}_p \rightarrow {\cal S}_{p-1}$ is defined by
\begin{eqnarray*}
L(x)(y) \equiv L(x,y) &:=& \frac{1}{2} \sum\limits^d_{i,j=1} (\sigma
\sigma^t)_{ij}(x,y) \partial^2_{ij}y
 -\sum\limits^d_{j=1} b_i (x,y) \partial_i y \\
\end{eqnarray*}

where $\sigma_{ij},b_i: \mathbb R^d \times {\cal S}_p \rightarrow
\mathbb R, i,j = 1,\cdots,d$ are assumed to satisfy a Lipschitz
condition as follows : Let $f: \mathbb R^d \times {\cal S}^\prime
\rightarrow \mathbb R$. We say that $f$ satisfies a $(p,q)$ local
Lipschitz condition , uniformly in $x \in \mathbb R^d$ if for all
$\lambda
>0$ there exists $C=C(\lambda,p,q)$ such that
\[
\left| f(x,\varphi) - f(x,\psi)\right|  \leq C\|\varphi -\psi\|_q
\]
for all $\varphi,\psi \in B_p (0,\lambda)$ and $x \in \mathbb R^d$.

Under the above condition, we can show the existence and uniqueness
of solutions of the above equation (\cite{BRASV}). Here, given a
measurable map $u : \mathbb R^d \rightarrow {\cal S}_p$ we will
assume the existence of a unique solution to the above PDE i.e. for
each $x \in \mathbb R^d$, the existence of a unique  map $u(.,x) :
[0,T] \rightarrow {\cal S}_p$ which is continuous and satisfies
\[ u(t,x) = u(x)+ \int\limits_0^t
L(x,u(s,x))ds.
\]
where the equation holds in ${\cal S}_q, q \leq p-1$. Suppose now we
are given a potential function i.e. a real valued function of the
form $c(t,x,y), 0 \leq t \leq T, x \in \mathbb R^d, y \in
C([0,T],{\cal S}_p)$, satisfying for each $x$, conditions 1 and 2 of
Section 2 for $H = {\cal S}_p$. Let $\hat \sigma_{ij}(t,x,.), \hat
b_i(t,x,.), \hat c(t,x,.)$ be as defined in \\ Section 3. For $t \in
[0,T], x \in \mathbb R^d, y \in C([0,T],{\cal S}_p)$ define the
operator
\begin{eqnarray*}
\hat L(s,x,y) &:=& \frac{1}{2} \sum\limits^d_{i,j=1} (\hat \sigma
\hat \sigma^t)_{ij}(s,x,y) \partial^2_{ij} y_s
 -\sum\limits^d_{j=1} \hat b_i (s,x,y) \partial_i y_s \\ &&+ \hat
 c(s,x,y)y_s.
\end{eqnarray*} The following theorem can be proved in the same
manner as Theorem (3.1).

\begin{Theorem} Let $(u(t,x))$ be a solution of equation (4.6) for a given
$u : \mathbb R^d \rightarrow {\cal S}_p$ and let $c(.,.)$ satisfy
conditions 1) and 2). Then,
\[
\hat u(t,x):= u(t,x)e^{ \int\limits_0^t c(s,x,u(.,x))ds}\] satisfies

\begin{eqnarray}
\partial_t \hat u(t,x) &=& \hat L(t,x,u(.,x))\\
\hat u(0,x) &=& u(x) \nonumber .
\end{eqnarray}
If equation (4.6) has a unique solution so has equation (4.7).
\end{Theorem}

For a given r-dimensional Brownian motion $(B_t)$ and $u(t,x)$
satisfying (4.6) let
\[
Z_t^x := \int\limits_0^t \sigma (x,u(s,x)) \cdot dB_s
+\int\limits_0^t b(x,u(s,x)) ds
\]
Let $Y_t^x := \tau_{Z_t^x}u(x)$, where $\tau_x : {\cal S}_p
\rightarrow {\cal S}_p$ are the translation operators. Note that for
each $x, E\|\tau_{Z_t^x}u(x)\|_p < \infty$. Then as in the proof of
Theorem 6.3, \cite{BR2}, we have $u(t,x) = EY_t^x$. Let $\hat Y_t^x
:= (\tau_{Z_t^x}u(x))e^{ \int\limits_0^t c(s,x,u(.,x))ds} $. Then we
have the following
\begin{Corollary} For each $x \in \mathbb R^d$, we have $$\hat u(t,x)
= E(Y_t^x)e^{ \int\limits_0^t c(s,x,EY^x_\cdot))ds}, 0 \leq t \leq
T. $$

\end{Corollary}

\begin{Remark} We note that for fixed $(t,x), \hat u(t,x) \in {\cal S}_p$
whenever $u(t,x) \in {\cal S}_p$. Further, Theorem (4.1) implies
that the degree of smoothness of $\hat u(t,x)$ in the backward
variable $x$ is the minimum of the degree of smoothness of the maps
$x \rightarrow u(t,x)$ and that of $x \rightarrow c(t,x,u(t,x))$.
\end{Remark}

\section{Conclusion}
In this section we make a few remarks on the applications of Theorem
2.2. We first consider the PDE (4.6) and its interplay with the
${\cal S}_p$ valued processes considered in Section 3. The existence
and uniqueness of solutions of (4.6) in the non-linear case will be
considered in a separate paper (\cite{BRASV}). Here we will consider
two separate classes of equation (4.6), in remarks 1 and 2 below,
corresponding to different classes of coefficients $\sigma_{ij},
b_i$ and the corresponding classes of linear operators $L(x,\phi)$
in (4.6).
\begin{enumerate}
\item{} We assume that the coefficients depend only on $x \in \mathbb R^d$
i.e $\sigma_{ij}(x,\phi) = \sigma_{ij}(x), b_i(x,\phi) = b_i(x), i=
1, \cdots, d, j = 1, \cdots,r$ and the initial condition $u(x)$ is
arbitrary. In this case $L(x) : {\cal S}_p \rightarrow {\cal S}_q, q
\leq p-1 $ is a linear operator. The solution $u(t,x)$ exists
uniquely - because of the monotonicity inequality satisfied by
$L(x)$ and is given by $u(t,x):= E\tau_{Z_t^x}u(x)$ where for each
$x \in \mathbb R^d$,
\[
Z_t^x :=  \sigma (x)\cdot B_t +  b(x)t
\]
In particular, $\hat u(t,x)$ is the unique solution to (4.7) for any
given potential function $c(t,x,y)$ satisfying conditions 1 and 2.
In this example, the role of the variable $x \in \mathbb R^d$ in the
coefficients of the equation is that of an `external parameter' and
as a consequence $(Z_t^x)$ is a Gaussian process, for each $x$.

\item{}
Suppose that $\bar \sigma_{ij}, \bar b_i \in {\cal S}_{-p}, p >
\frac{d}{4}$ and $\sigma_{ij}, b_i : {\cal S}_p \rightarrow \mathbb
R$ are defined by $\sigma_{ij}(\phi) := \langle \bar \sigma_{ij},
\phi \rangle$ etc. In particular $\bar \sigma_{ij}, \bar b_i $ do
not depend on $\phi \in {\cal S}_p$ and $\sigma_{ij}, b_i  $ do not
depend on $x \in \mathbb R^d$. Consider the operators $L,\bar L$,
respectively non-linear and linear, associated,respectively with
(3.4) and (4.6). In the following computations we will show the
connection between solutions of (3.4)-(3.5) associated with the
non-linear operator $L$ and the solutions of (4.6)-(4.7) associated
with the linear operator $L(x)$ which we here denote by $\bar L(x)$,
acting on $\phi \in {\cal S}_p, p > \frac{d}{4} $ as follows :

\[\bar L\phi(x) := \frac{1}{2} \sum\limits^d_{i,j=1} (\bar \sigma
 \bar \sigma^t)_{ij}(x) \partial^2_{ij} \phi(x)
 +\sum\limits^d_{j=1} \bar b_i (x) \partial_i \phi(x).\]

If $\sigma_{ij},b_i$ are bounded measurable functions and $p
> 1$ then $\bar L : {\cal S}_p \rightarrow {\cal S}_0 = L^2$.
Associated (as above) with the coefficients $\bar \sigma_{ij}, \bar
b_i \in {\cal S}_p$ is the non-linear operator $L : {\cal S}_p
\rightarrow {\cal S}_{p-1}$

\[L(\phi) := \frac{1}{2} \sum\limits^d_{i,j=1} ( \sigma
 \sigma^t)_{ij}(\phi) \partial^2_{ij} \phi
 +\sum\limits^d_{j=1} b_i (\phi) \partial_i \phi.\]

where $\sigma_{ij}, b_i : {\cal S}_p \rightarrow \mathbb R$ are
defined above. Let $(X_t^x)$ denote the unique solution of
\begin{eqnarray*}
dX_t &=& \bar \sigma (X_t) \cdot dB_t +\bar b(X_t) dt \\
X_0 &=& x.
\end{eqnarray*}
Let $Y_t := \delta_{X_t^x} \in {\cal S}_{-p}, t \geq 0,$ and for
simplicity we assume the associated life time $\eta^x = \infty$
almost surely. We denote by $(P_t)$ the transition semi-group
corresponding to $(X_t^x)$ and by $P_t^\ast(x)$ the kernel
$P_t^\ast(x) := EY_t^x = E\delta_{X_t^x}$ which is just the
transition probability measure of $(X_t^x)$ represented as an
element of ${\cal S}_{-p}$. Taking expected values in (3.4) we see
that $P_t^\ast(x)$ satisfies
\begin{eqnarray}
\partial_t u(t,x) &=& \bar L^\ast u(t,x) \\
u(0,x) &=& \delta_x. \nonumber
\end{eqnarray}
 where $\bar L^\ast(x)$ is the formal adjoint of $\bar L$ satisfying :
 \[  \bar L^\ast P_t^\ast (x) = EL(\delta_{X_t^x}).\] This maybe
verified by acting with a test function $u \in {\cal S}$. Equation
(5.8) is the same as equation (4.6) for the linear operator $L(x) =
\bar L^\ast(x)$. When $\sigma_{ij}, b_i$ are twice continuously
differentiable with bounded derivatives then (5.8) has a unique
solution (see \cite{DS}, Theorem 2.2.9).On the other hand, the
operator $\hat L$ in (4.7) when $L(x)= \bar L^\ast(x)$ is just $\bar
L^\ast(x) + V^\ast(x) = \bar L^\ast + V(x)$ and hence the solution
of (4.6) viz. $(P_t^\ast(x))$ transforms into the solution of (4.7)
viz. $P_t^\ast(x)e^{tV(x)}$. Hence by the uniqueness result quoted
above and the uniqueness result in Theorem (4.1), the evolution
equation
\begin{eqnarray}
\partial_t u(t,x) &=& (\bar L^\ast+V) u(t,x) \\
u(0,x) &=& \delta_x. \nonumber
\end{eqnarray}
has a unique solution given by $\hat u(t,x):=
P_t^\ast(x)e^{tV(x)}$,when $\sigma_{ij}, b_i$ are twice continuously
differentiable with bounded derivatives .

For $V \in {\cal S}_p$ define $c(.,.) : [0,\infty) \times
C([0,\infty),{\cal S}_{-p}) \rightarrow \mathbb R$ as $c(t,y) :=~\\
<V,y(t)> ,~ y \in C([0,\infty),{\cal S}_{-p}).$ Let
\[
\hat Y_t := Y_t~e^{\int\limits_0^t c(s,Y)ds}
\]
where $c(s,Y):=\langle V,Y_s\rangle=V(X^x_s)$. Since  $|c(t,Y)| \leq
\|V\|_p K$ where $\|\delta_z\|_{-p} \leq K$, we have
\[
E\|\hat Y_t\|_p \leq e^{K\|V\|_pt} E\|Y_t\|_p < \infty.
\]

Let for $f \in {\cal S}_p$, \[ P^{V}_tf(x) :=E(e^{\int\limits_0^t
V(X^x_s)ds}f(X_t^x))
\] and let \[ P^{V\ast}_t(x) := E\hat Y_t =E(e^{\int\limits_0^t
V(X^x_s)ds}\delta_{X_t^x}) \in {\cal S}_{-p}. \]

Then the following calculations show that $P^{V\ast}_t(x)$ satisfies
(5.9).  Let $f \in {\cal S}$. From the definition of $\hat L(t,y)$
we have
\begin{eqnarray*}
\langle f, \hat L(t,\hat Y)\rangle &=& e^{\int\limits_0^t
V(X_s^x)ds}\langle f, L(\delta_{X_t^x})\rangle
 + e^{\int\limits_0^t V(X_s^x)ds}V(X_t^x) \langle f,\delta_{X_t^x}\rangle\\
&=& e^{\int\limits_0^t V(X_s^x)ds}(\bar L+V) f(X_t^x).
\end{eqnarray*}
Hence from the equation satisfied by $\hat Y_t$ we get
\begin{eqnarray*}
 \langle f, P^{V\ast}_t(x)\rangle &=& P_t^V f(x)= \langle f,E\hat
Y_t\rangle
= f(x) + \int_0^t \langle f, \hat L(s,\hat Y)\rangle ds \\
&=& f(x) +\int\limits_0^t E[e^{\int\limits_0^s V(X_u^x)du}(\bar L+V) f(X^x_s)] ds \\
&=&  f(x) +\int\limits_0^t P_s^V((\bar
L+V)f)(x) ds \\
&=& \langle f, \delta_x\rangle  +\int\limits_0^t \langle f, \bar
L^\ast+V)(P^{V\ast}_s(x) \rangle ds.
\end{eqnarray*}
It follows by uniqueness of solutions of (5.9) that with the
coefficients $\sigma_{ij}, b_i$ as above, we have the following
special case of Corollary (4.2) with $c(t,x,y):= V(x)$ and with
equality in ${\cal S}_{-p}$, for each $x \in \mathbb R^d$ :
\[P^{V\ast}_t(x) = P^{\ast}_t(x)e^{tV(x)} =
e^{tV(x)}E\delta_{X_t^x}.\]

\end{enumerate}

}

\end{document}